\documentclass{amsart}
\usepackage[centertags]{amsmath}
\usepackage{amsfonts}
\usepackage{amssymb}
\usepackage{amsthm}
\usepackage{amsmath}
\usepackage{amscd}
\usepackage{latexsym}
\usepackage{color}
\usepackage{graphics}

\theoremstyle{plain}
\newtheorem{thm}{Theorem}[section]
\newtheorem{cor}[thm]{Corollary}
\newtheorem{lem}[thm]{Lemma}

\theoremstyle{definition} \theoremstyle{remark}
\numberwithin{equation}{section}

\newcommand{\ep}{\varepsilon}

\newcommand{\vph}{\varphi}

\newcommand{\intl}{\int\limits}

\begin{document}

%\selectlanguage{english}\fontencoding{T2A}\selectfont

\title[discrete magnetic Sch\"{o}dinger operator]{Essential self-adjointness of a discrete magnetic Sch\"{o}dinger operator}
\author{Volodymyr Sushch}
%\shorttitle{Short paper title for the headers}
%\shortauthor{F. Author, S. Author}
%\date{28.02.2005}
\address{Department of Mathematics, Koszalin University of Technology, Sniadeckich 2,
 75-453 Koszalin, Poland;  Pidstrygach Institute for Applied Problems of Mechanics and
 Mathematics NASU,  Lviv, Ukraine }

\email{volodymyr.sushch@tu.koszalin.pl}
 \maketitle
\begin{abstract}
 We prove essential self-adjointness for a semibounded from below discrete magnetic Schr\"{o}dinger operator
in a space that represents a combinatorial model of the
two-dimensional Euclidean space. The Dezin discretization scheme is
used for constructing a discrete model.
\end{abstract}

\section{Introduction}
  Let $\Lambda^p_k(\mathbb{R}^2)$ be the set of all $k$-smooth (i.e., of the class  $C^k$) complex-valued $p$-forms in
$\mathbb{R}^2$ and let
$\Lambda^p(\mathbb{R}^2)=\Lambda^p_\infty(\mathbb{R}^2)$.
We define a magnetic potential as a real-valued 1-form  $A\in
\Lambda^1_1(\mathbb{R}^2)$, i.e.,
$$A=A_1dx^1+A_2dx^2,$$ where $A_1, A_2\in C^1(\mathbb{R}^2)$  are
real-valued functions. We introduce an invariant inner product for
$p$-forms with compact support in $\mathbb{R}^2$ in the following
way

\begin{equation} \label{1}(\vph, \
\psi)=\intl_{\mathbb{R}^2}\vph\wedge\ast\overline{\psi},
\end{equation}
where $\ast$ is the operation of metric conjugation of forms (the
Hodge star operator), $\wedge$ is the operation of exterior
multiplication and the bar over $\psi$ denotes complex conjugation.
Consider the completion of the linear spaces of smooth forms in norm
that is generated by the inner product (\ref{1}). We denote the
formed Hilbert spaces by
 $L^2(\mathbb{R}^2)$ for 0-forms (functions) and by
$L^2\Lambda^p(\mathbb{R}^2)$ for $p$-forms, $p = 1,2$. Let $d$ be
the operator of exterior differentiation. We introduce a deformed
differential according to the rule
\begin{equation}\label{2}
 d_A:C^{\infty}(\mathbb{R}^2)\rightarrow
\Lambda^{1}_1(\mathbb{R}^2), \quad \vph\rightarrow d\vph+i\vph A,
 \end{equation}
where $i^2=-1$ and $A$ is the magnetic potential. The inner product
(\ref{1}) enables us to define an operator formally adjoint to
$d_A$ as follows
$$\delta_A:\Lambda^{1}_1(\mathbb{R}^2)\rightarrow
C(\mathbb{R}^2).$$ Then we may define the magnetic Laplacian
$\Delta_A$ (Laplacian $\Delta$ with potential A) according to
\begin{equation}\label{3}
-\Delta_A\equiv\delta_Ad_A:C^{\infty}(\mathbb{R}^2)\rightarrow
C(\mathbb{R}^2).
 \end{equation}
 Identifying the magnetic potential $A$ with the multiplication operator
\begin{equation*}
 A: C^{\infty}(\mathbb{R}^2)\rightarrow
\Lambda^{1}_1(\mathbb{R}^2), \qquad \vph\rightarrow\vph A,
\end{equation*}
we may represent the operator  $\delta_A$ in the form
\begin{equation}
\label{3.1} \delta_A\omega=(\delta-iA^\ast)\omega,
\end{equation}
where $\delta, \ A^\ast$ are operators formally adjoint to $d$ and $A$, respectively. Using (\ref{2}) and (\ref{3.1}), we may rewrite
the magnetic Laplacian $\Delta_A$ as
 \begin{align*}
-\Delta_A\vph &=(\delta-iA^\ast)(d\vph+iA\vph)=\\
&=-\Delta\vph-iA^\ast d\vph+i\delta(A\vph)+A^\ast A\vph.
\end{align*}
Consider now the magnetic Schr\"{o}dinger operator
\begin{equation}\label{4}
H_{A,V}=-\Delta_A+V,
\end{equation}
where $V$ is a real-valued function, which is also called electric potential, and $V\in
L^2_{loc}(\Bbb R^2)$.
 Suppose that the operator  $H_{A,V}$ is semi-bounded from below on $C^{\infty}_{0}(\mathbb{R}^2)$,
 i.e., there exists a constant $c\in\mathbb{R}$ such that

\begin{equation}\label{5}
(H_{A,V}\varphi, \ \varphi)\geq-c(\varphi, \ \varphi), \quad
\varphi\in C^{\infty}_{0}(\mathbb{R}^2).
\end{equation}
Here $C^{\infty}_{0}(\mathbb{R}^2)$ is the set of all  $C^{\infty}$ functions with compact support in $\mathbb{R}^2$. Then, as is well known (see
 \cite{Shubin}), operator (\ref{4}) is essentially self-adjoint on $C^{\infty}_{0}(\mathbb{R}^2)$.

The main aim of the present work is to study the essential self-adjointness of the discrete magnetic
Schr\"{o}dinger operator on a combinatorial object corresponding to
$C^{\infty}_{0}(\mathbb{R}^2)$. In \cite{Sushch}, we proposed a discrete
model of the magnetic Laplacian (\ref{3}) such that it preserves the geometrical structure of the initial continual object
and also proved the self-adjointness of the operator of a discrete Dirichlet problem for the magnetic Laplacian.
In bounded domains, which gives the finite dimensionality of the corresponding Hilbert spaces of a discrete
problem, the results of \cite{Sushch} can easily be generalized for the case of the magnetic Schr\"{o}dinger operator.
In the present work we show that the semi-bounded from below discrete magnetic Schr\"{o}dinger operator, as
in the continual case, has a unique self-adjoint realization. It should be emphasized that, besides conditions (\ref{5}),
no other limitations are imposed on the behavior of discrete analogs of the potentials $A$ and $V$ at infinity.
Our approach is based on the formalism proposed by Dezin in \cite{Dezin}. We shall also use the results described in \cite{Sushch}.

Note that the discrete magnetic Laplacian and discrete magnetic Schr\"{o}dinger operators are fairly popular
subjects of inquiry among mathematicians and physicists. There exist numerous various approaches (different
from that proposed in the present work) both to the construction of discrete models and to investigation of the
corresponding difference operators (see, e.g., \cite{BSBE, Pastur, DodM, Higuchi,
Shubin1, Sunada} and references therein). In the overwhelming majority of
these works particular attention is given to studying the spectral properties of discrete magnetic Schr\"{o}dinger
operators on infinite graphs. As to investigations of the essential self-adjointness of discrete operators a review
of different aspects of this problem can be found in \cite{Berez, Teschl}.

\section{Basic Combinatorial Constructions}

In this section we briefly recall the definitions of the basic combinatorial operations, which will be used in
constructing the discrete analogs of operators (\ref{3}) and (\ref{4}). Let
$\mathfrak{C}(2)$ be a two-dimensional complex, i.e., a
combinatorial model of $\mathbb{R}^2$ (for more detail, see \cite{Dezin}, \cite{Sushch}). The complex
$\mathfrak{C}(2)$ can be represented as
$\mathfrak{C}(2)=\mathfrak{C}^0\oplus\mathfrak{C}^1\oplus\mathfrak{C}^2$, where  $\mathfrak{C}^p$ is a real linear space of $p$-dimensional chains,  $p=0,1,2$. We denote by
$\{x_{k,s}\}, \ \{e_{k,s}^1, \
e_{k,s}^2\}, \ \{\Omega_{k,s}\},$ \ $k, s\in\mathbb{Z},$ the sets of basic elements of the spaces $\mathfrak{C}^0, \ \mathfrak{C}^1$  and
$\mathfrak{C}^2$, respectively.
For convenience, we introduce the shift operators
 $$\tau k=k+1, \qquad \sigma k=k-1$$
on the set of indices. The boundary operator $\partial$ on the basic elements of
$\mathfrak{C}(2)$ is assigned as
\begin{align}\notag
\partial x_{k,s}=0, &\qquad \partial e_{k,s}^1=x_{\tau k,s}-x_{k,s}, \qquad
\partial e_{k,s}^2=x_{k,\tau s}-x_{k,s}, \\
\label{6} &\partial\Omega_{k,s}=e_{k,s}^1+e_{\tau k,s}^2-e_{k,\tau
s}^1-e_{k,s}^2.
\end{align}

We introduce a dual object, i.e., a complex conjugate to
$\mathfrak{C}(2)$. We denote it by $K(2)$, and let it be a linear
space of complex-valued functions over
$\mathfrak{C}(2)$. Further, suppose that $K^0, \  K^1$ and
$K^2$ are linear spaces
conjugate to $\mathfrak{C}^0,
\mathfrak{C}^1$ and $\mathfrak{C}^2$,  i.e., they have bases of the form $\{x^{k,s}\}, \ \{e^{k,s}_1, \ e^{k,s}_2\}, \ \{\Omega^{k,s}\}$, respectively.
Then we may consider $K(2)=K^0\oplus K^1\oplus K^2$ as a complex of complex-valued cochains of the corresponding
dimensionality. In what follows, these cochains are called forms, which emphasizes their proximity to
the corresponding continual objects (differential forms). Then the 0-, 1-, and 2-forms $\vph\in K^0, \
\omega=(u, v)\in K^1$ and  $\eta\in K^2$ look like
\begin{equation}\label{7}
\vph=\sum_{k,s}\vph_{k,s}x^{k,s}, \quad
\omega=\sum_{k,s}(u_{k,s}e^{k,s}_1+v_{k,s}e^{k,s}_2), \quad
\eta=\sum_{k,s}\eta_{k,s}\Omega^{k,s},
\end{equation}
where $\vph_{k,s}, \ u_{k,s}, \
v_{k,s},  \ \eta_{k,s}\in\mathbb{C}$ for all $k,s\in\mathbb{Z}$.

We define the operation of pairing for the basic elements of complexes $\mathfrak{C}(2)$ and $K(2)$ according to the
rule
\begin{equation}\label{8}
<x_{k,s}, \ x^{p,q}>=<e_{k,s}^1, \ e^{p,q}_1>=<e_{k,s}^2, \
e^{p,q}_2>=<\Omega_{k,s}, \ \Omega^{p,q}>=\delta_{k,s}^{p,q},
\end{equation}
where $\delta_{k,s}^{p,q}$ is the Kronecker delta. Pairing (\ref{8}) is extended to arbitrary forms (\ref{7}) by linearity.
The boundary operator (\ref{6}) induces in the conjugate complex $K(2)$  a dual operation, namely, a coboundary
operator $d^c$:
\begin{equation}\label{9}
<\partial a, \ \alpha>=<a, \ d^c\alpha>,
\end{equation}
where $a\in \mathfrak{C}(2)$ and  $\alpha\in K(2)$. We consider the coboundary operator 
\begin{equation*}d^c: K^p \rightarrow K^{p+1}
\end{equation*}
 as a discrete analog
of the operation of exterior differentiation $d$. In what follows, we use the following difference representations
of the operator $d^c$:
\begin{align}\notag
 &<e^1_{k,s}, \ d^c\vph>=\vph_{\tau k,
 s}-\vph_{k,s}\equiv\Delta_k\vph_{k,s},\\
 \label{10}
 & <e^2_{k,s}, \ d^c\vph>=\vph_{k, \tau
 s}-\vph_{k,s}\equiv\Delta_s\vph_{k,s},\\ \notag
  &<\Omega_{k,s}, \ d^c\omega>=v_{\tau k,
 s}-v_{k,s}-u_{k,\tau s}+u_{k,s}\equiv\Delta_k v_{k,s}-\Delta_s
 u_{k,s}.
 \end{align}

In the complex $K(2)$ we introduce the operation of multiplication, which is considered as an analog of
the exterior multiplication of differential forms. We denote this operation by $\cup$ and define it as
\begin{align}\label{11}\notag
&x^{k,s}\cup x^{k,s}=x^{k,s}, \qquad e^{k,s}_2\cup e^{k,\tau
s}_1=-\Omega^{k,s},\\  &x^{k,s}\cup e^{k,s}_1=e^{k,s}_1\cup x^{\tau
k,s}=e^{k,s}_1, \qquad x^{k,s}\cup e^{k,s}_2=e^{k,s}_2\cup x^{k,\tau
s}=e^{k,s}_2,
\\ \notag &x^{k,s}\cup \Omega^{k,s}= \Omega^{k,s}\cup x^{\tau k,\tau s}=e^{k,s}_1\cup
e^{\tau k,s}_2=\Omega^{k,s},
 \end{align}
assuming that the product is equal to zero in all other cases.
The $\cup$-multiplication is extended to forms (\ref{7}) by linearity.
We denote by $\ep^{k,s}$ an arbitrary basic element of $K(2)$. Then we introduce an operation $\ast$ taking
\begin{equation}\label{12}
\ep^{k,s}\cup\ast\ep^{k,s}=\Omega^{k,s}.
\end{equation}
Using (\ref{11}), we obtain
\begin{equation*}
\ast x^{k,s}=\Omega^{k,s},  \quad \ast e^{k,s}_1=e^{\tau k,s}_2,
\quad \ast e^{k,s}_2=-e^{k,\tau s}_1,  \quad \ast
\Omega^{k,s}=x^{\tau k,\tau s}.
\end{equation*}
The operation $\ast$ is extended to arbitrary forms by linearity.

Let $\alpha\in K^p$ be an arbitrary $p$-form, i.e.,
\begin{equation}\label{13}
\alpha=\sum_{k,s}\alpha_{k,s}\ep^{k,s}.
\end{equation}
We denote by $K^p_0$
 the set of all discrete $p$-forms with compact support, i.e., if $\alpha\in K^p_0$, then only a finite
number of components $\alpha_{k,s}$  in  (\ref{13})  is different from zero. Let now
\begin{equation}\label{14}
\Omega=\sum_{k,s}\Omega_{k,s},  \qquad k,s\in\mathbb{Z},
\end{equation}
where $\Omega_{k,s}$ is a two-dimensional basic element of $\mathfrak{C}(2)$. Note that we shall also use the notation $\Omega=\Omega_N$, if
sum (\ref{14}) is finite and $-N\leq k,s\leq N,$   $N\in\mathbb{N}$.

The relation
\begin{equation}\label{15}
(\alpha, \ \beta)=<\Omega, \ \alpha\cup\ast\overline{\beta}>,
\end{equation}
where $\alpha, \beta\in K^p_0$, gives a correct definition of the inner product for discrete $p$-forms (see (\ref{1})). Using (\ref{8}), (\ref{11})
 and (\ref{12}), we may rewrite relation (\ref{15}) as
\begin{equation}\label{16}
(\alpha, \ \beta)=\sum_{k,s}\alpha_{k,s}\overline{\beta_{k,s}}.
\end{equation}

The inner product (\ref{15}) enables us to define an operator formally adjoint to $d^c$, namely, the operator
$\delta^c:
K^{p+1}\rightarrow K^p$  satisfying the following relation
\begin{equation}\label{17}
(d^c\alpha, \ \beta)=(\alpha, \ \delta^c\beta), \quad \alpha\in
K^p_0, \quad \beta\in K^{p+1}.
\end{equation}
It is easy to show that
\begin{equation}\label{18}
\delta^c\beta=(-1)^{p+1}\ast^{-1}d^c\ast\beta,
\end{equation}
where $\ast^{-1}$ is the operation inverse to $\ast$, i.e., $\ast^{-1}\ast=1$.
 Hence, we may consider the operator $\delta^c$  as a
discrete analog of the codifferential $\delta$. Taking (\ref{10}) into account, we have for $\omega\in K^1$
\begin{equation}\label{19}
\delta^c\omega=\sum_{k,s} (-\Delta_k u_{\sigma k,s}-\Delta_s
v_{k,\sigma s})x^{k,s}.
\end{equation}
Thus, the discrete analog of the Laplace operator looks like

\begin{equation}\label{20}
-\Delta^c=\delta^c d^c+d^c\delta^c: K^p\rightarrow K^p.
\end{equation}
Obviously, since $\delta^c \vph=0$ for $\vph\in K^0$, we have
\begin{equation}\label{21}
-\Delta^c\vph=\delta^c d^c\vph.
\end{equation}

\section{Discrete Analog of the Magnetic Laplacian}
Let a real-valued 1-form
\begin{equation*}
 A=\sum_{k,s}(A_{k,s}^{1}e_{1}^{k,s}+A_{k,s}^{2}e_{2}^{k,s}),
\end{equation*}
where
$A_{k,s}^{1}, A_{k,s}^{2}\in\mathbb{R},$ be a discrete analog of the magnetic potential. Then we define a discrete analog of the
deformed differential (\ref{2}) in the following way:
\begin{equation}\label{22}
d^c_A: K^0\rightarrow K^1, \qquad \vph\rightarrow d^c\vph+i\vph\cup
A.
\end{equation}
In view of (\ref{10}) and (\ref{11}), we obtain
\begin{equation}\label{23}
d^c_A\vph=\sum_{k,s}\big((\Delta_k\vph_{k,s}+i\vph_{k,s}A_{k,s}^{1})e_{1}^{k,s}+
(\Delta_s\vph_{k,s}+i\vph_{k,s}A_{k,s}^{2})e_{2}^{k,s}\big).
\end{equation}
Further, we identify the discrete magnetic potential $A$  with the operator of multiplication as follows:
\begin{equation}\label{24}
A:K^0\rightarrow K^1,  \qquad \vph\rightarrow\vph\cup A.
\end{equation}
Then it is easy to obtain
\begin{equation*}
 A\vph=\sum_{k,s}(\vph_{k,s}A_{k,s}^{1}e_{1}^{k,s}+\vph_{k,s}A_{k,s}^{2}e_{2}^{k,s}).
\end{equation*}
Let $A^\ast:K^1\rightarrow K^0$
 be the operator formally conjugate to $A$, i.e., it acts on an arbitrary 1-form $\omega=(u, v)$
according to the rule
\begin{equation}\label{25}
A^{\ast}\omega=\sum\limits_{k,s}(A^{1}_{k,s}u_{k,s}+A^{2}_{k,s}v_{k,s})x^{k,s}.
\end{equation}
Hence (see \cite{Sushch} for more detail), the operator $\delta^c_A:K^1\rightarrow K^0$,
 which is formally adjoint to the operator $d^c_A$, has
the form
\begin{equation}\label{26}
\delta^c_A=\delta^c-iA^\ast.
\end{equation}
Thus, we may define a discrete magnetic Laplacian as
\begin{equation*}
-\Delta^c_A=\delta^c_A d^c_A: K^0\rightarrow K^0.
\end{equation*}
In view of (\ref{22}) and (\ref{26}), we obtain
\begin{align}\label{27}\notag
-\Delta^c_A\vph&=\delta^c_A(d^c\vph+i\vph\cup A)=\\\notag
&=(\delta^c-iA^\ast)d^c\vph+(\delta^c-iA^\ast)(i\vph\cup A)=\\\notag
&=-\Delta^c\vph-iA^\ast d^c\vph+i\delta^c(\vph\cup A) +
A^\ast(\vph\cup A)=\\ &=-\Delta^c\vph-iA^\ast d^c\vph+i\delta^cA\vph
+ A^\ast A\vph.
\end{align}
Using (\ref{11}) and (\ref{19}), it is easy to show that for forms $\varphi\in K^0$ and
$\omega\in K^1$ (see (\ref{7})) one may write
\begin{equation*}
\delta^c(\omega\cup\varphi)=\delta^c\omega\cup\varphi-\sum_{k,s}\big(u_{\sigma
k,s}(\Delta_k\varphi_{k,s})+v_{k, \sigma
s}(\Delta_s\varphi_{k,s})\big)x^{k,s},
\end{equation*}
\begin{equation*}
\delta^c(\varphi\cup\omega)=\varphi\cup\delta^c\omega-\sum_{k,s}\big((\Delta_k\varphi_{\sigma
k,s})u_{\sigma k,s}+(\Delta_s\varphi_{k,\sigma s})v_{k, \sigma
s}\big)x^{k,s}.
\end{equation*}
Then, in view of (\ref{25}) and (\ref{26}), the discrete analogs of the Leibniz rule for the operator  $а\delta_A^c$
 will have the form
\begin{align}\label{28}\notag
\delta_A^c(\omega\cup\varphi)=\delta^c\omega\cup\varphi-&\sum_{k,s}\big(u_{\sigma
k,s}(\Delta_k\varphi_{k,s})+v_{k, \sigma
s}(\Delta_s\varphi_{k,s})\big)x^{k,s}-\\
-&i\sum_{k,s}\big(A_{k,s}^1u_{k,s}\varphi_{\tau
k,s}+A_{k,s}^2v_{k,s}\varphi_{k,\tau s}\big)x^{k,s},
\end{align}
\begin{equation}\label{29}
\delta_A^c(\varphi\cup\omega)=\varphi\cup\delta_A^c\omega-\sum_{k,s}\big((\Delta_k\varphi_{\sigma
k,s})u_{\sigma k,s}+(\Delta_s\varphi_{k,\sigma s})v_{k, \sigma
s}\big)x^{k,s}.
\end{equation}
In addition, we have for $\varphi,\ \psi\in K^0$
\begin{align*}
-\Delta_A^c(\varphi\cup\psi)&=\delta_A^c(d^c(\varphi\cup\psi)+i\varphi\cup\psi\cup
A)=\\&=\delta_A^c(d^c\varphi\cup\psi+\varphi\cup d^c\psi+
i\varphi\cup(\psi\cup
A))=\\&=\delta_A^c(d^c\varphi\cup\psi)+\delta_A^c(\varphi\cup
d_A^c\psi).
\end{align*}
From here, replacing $\omega$ by the 1-form $d^c\varphi$ (see (\ref{10})) and $\varphi$ by $\psi$ in (\ref{28}) as well as $\omega$ by the 1-form
$d_A^c\psi$, looking like (\ref{23}), in (\ref{29}), we obtain
\begin{equation}\label{30}
-\Delta_A^c(\varphi\cup\psi)=\varphi\cup\delta_A^cd_A^c\psi+\delta^cd^c\varphi\cup\psi+\sum_{k,s}\Phi_{k,s}x^{k,s},
\end{equation}
where
\begin{align}\label{31}\notag
\Phi_{k,s}&=(\Delta_k\varphi_{\sigma k,s})(\psi_{\tau
k,s}-\psi_{\sigma k,s}+i\psi_{\sigma k,s}A^1_{\sigma
k,s})+i(\Delta_k\varphi_{k,s})\psi_{\sigma k,s}A^1_{k,s}+\\&+
(\Delta_s\varphi_{k,\sigma s})(\psi_{k,\tau s}-\psi_{k,\sigma
s}+i\psi_{k,\sigma s}A^2_{k,\sigma
s})+i(\Delta_k\varphi_{k,s})\psi_{k,\sigma s}A^2_{k,s}.
\end{align}

\section{Discrete Magnetic Schr\"{o}dinger Operator}

Let a real-valued 0-form $V\in K^0$ be a discrete analog of the
electric potential, i.e.,
 \begin{equation*}
 V=\sum_{k,s}V_{k,s}x^{k,s}, \qquad V_{k,s}\in\mathbb{R}.
\end{equation*}
Then the discrete analog of the magnetic Schrцdinger operator
(\ref{4}) has the form
 \begin{equation}\label{32}
 H_{A,V}^c=-\Delta^c_A+V.
\end{equation}
 Since we do not
impose any restrictions on the behavior of components of the
discrete forms $A$ and $V$ at infinity, operator (\ref{32}),
generally speaking, is unbounded.

For forms looking like (\ref{13}), we introduce a linear space
\begin{equation}\label{33} \mathcal{H}^p=\{\alpha\in K^p: \
\sum_{k,s}|\alpha_{k,s}|^2<+\infty, \quad k,s\in\mathbb{Z}\}, \qquad
p=0,1,2.
 \end{equation}

Obviously, according to (\ref{16}), the space $\mathcal{H}^p$ is a
Hilbert space with the inner product (\ref{15}) and a norm
\begin{equation}\label{34} \|\alpha\|=\sqrt{(\alpha, \
\alpha)}=\Big(\sum_{k,s}|\alpha_{k,s}|^2\Big)^{\frac{1}{2}}.
\end{equation}

It should be noted that, if $\alpha\in\mathcal{H}^0$, then the
sequence of components $(\alpha_{k,s})$ is an element of
$\ell^2(\mathbb{Z}^2)$, i.e., the space of all square summable
complex-valued sequences. Since the set of all finite sequences
$\ell_0(\mathbb{Z}^2)$ is dense in $\ell^2(\mathbb{Z}^2)$, the space
 $K_0^0$ is dense in  $\mathcal{H}^0$. Hence, the operator $H_{A,V}^c:
K_0^0\rightarrow\mathcal{H}^0$  is densely defined (i.e.,
$\overline{K_0^0}=\mathcal{H}^0$) and symmetric. In what follows,
we assume that operator (\ref{32}) is semi-bounded from below on
$K_0^0$ 0 , i.e., condition (\ref{5}) is satisfied for $H_{A,V}^c$ 
and all $\varphi\in K_0^0$ .

We define minimal and maximal operators associated with $H_{A,V}^c$
 in $\mathcal{H}^0$ as follows:
  \begin{equation*}
H_{min}: D(H_{min})\rightarrow\mathcal{H}^0, \qquad H_{max}:
D(H_{max})\rightarrow\mathcal{H}^0,
\end{equation*}
 where
 \begin{equation*}
D(H_{min})=K^0_0, \qquad D(H_{max})=\{\varphi\in\mathcal{H}^0| \
H_{A,V}^c\varphi\in\mathcal{H}^0\}.
\end{equation*}
The essential self-adjointness of the operator $H_{A,V}^c$ means
that $\overline{H_{min}}=H_{max}$, i.e., the closure of the minimal
operator in $\mathcal{H}^0$ coincides with the maximal operator.

Further, we introduce a cutting 0-form $\chi^N\in K^0$ by
    \begin{equation}\label{37}
\chi^N=\sum_{k,s}\chi^N_{k,s}x^{k,s},  \quad \mbox{де} \quad
\chi^N_{k,s}=\left\{\begin{array}{l}1,\quad |k|, |s|\leq N\\
                            0, \quad |k|, |s|>N
                            \end{array}\right., \quad N\in
                            \mathbb{N}.
\end{equation}

Hence, $\chi^N=\sum x^{k,s}$ and $k, s$ take values from $-N$ to
$N$. We also denote the inner product (\ref{15}) by $(\cdot ,
\cdot)_N$ if $\Omega=\Omega_N$ (see (\ref{14})).
 \begin{lem} Let
$\psi\in K^0$ and let $\chi^{\tau N}\in K^0$  be a form looking like
(\ref{37}). Then \begin{equation}\label{38}
\big(H_{A,V}^c(\chi^{\tau N}\cup\psi), \ \chi^{\tau
N}\cup\psi\big)_N=\big(\chi^{\tau N}\cup H_{A,V}^c\psi, \ \chi^{\tau
N}\cup\psi\big)_N.
\end{equation}
\end{lem}
\begin{proof}
  Using relation (\ref{30}), we obtain for arbitrary
$\varphi, \psi\in K^0$
\begin{align}\label{39}\notag
H_{A,V}^c(\varphi\cup\psi)&=\delta_A^cd_A^c(\varphi\cup\psi)+V\cup(\varphi\cup\psi)=\\
&=\varphi\cup
H_{A,V}^c\psi+\delta^cd^c\varphi\cup\psi+\sum_{k,s}\Phi_{k,s}x^{k,s}.
\end{align}
 Since the components of 
form \ $\delta^cd^c\varphi$  \ look like the following difference operators
$$-\Delta_k(\Delta_k\varphi_{\sigma
k,s})-\Delta_s(\Delta_s\varphi_{k,\sigma s})$$
 and all summands of
the components of $\Phi_{k,s}$ have multipliers of the form
$\Delta_k\varphi_{k,s}$ and  $\Delta_s\varphi_{k,s}$ (see (\ref{31})),
we find for $\varphi$ with constant components that
$\delta^cd^c\varphi=0$  and $\Phi_{k,s}=0$. Let now
$\varphi=\chi^{\tau N}$. We denote by $\Phi^{\tau N}$ the 0-form
with components (\ref{31}), where $\varphi_{k,s}$ are replaced by
$\chi^{\tau N}_{k,s}$. Substituting (\ref{39}) in the inner product
$\big(H_{A,V}^c(\chi^{\tau N}\cup\psi), \ \chi^{\tau
N}\cup\psi\big)_N$, we see that the components of the form $
\chi^{\tau N}$  are equal to 1 at points $x_{k,s}$
 of the domain $\Omega_N$ and by a step beyond its boundary. This fact
guarantees that the components of the forms $\delta^cd^c\chi^{\tau
N}$ and $\Phi^{\tau N}$ are equal to zero at points of the boundary
of the domain $\Omega_N$, i.e., for $k=\pm N$ or $s=\pm N$. From
here we immediately obtain
\begin{equation*}\big(\delta^cd^c\chi^{\tau N}\cup\psi, \ \chi^{\tau
N}\cup\psi\big)_N=0, \qquad \big(\Phi^{\tau N}, \ \chi^{\tau
N}\cup\psi\big)_N=0
\end{equation*} and this means that equality (\ref{39}) holds true.
\end{proof}

 \begin{thm}
 Let the discrete magnetic Schr\"{o}dinger operator $H_{A,V}^c$ be
semi-bounded from below on $K^0_0$. Then $H_{A,V}^c$ is essentially
self-adjoint.
\end{thm}
\begin{proof}
Obviously, every semi-bounded operator becomes strictly positive if
the corresponding constant is added to it. For example, adding
$(c+1)Id$ to $H_{A,V}^c$, we get
 \begin{equation*}
\big(H_{A,V}^c\psi, \ \psi\big)\geq\|\psi\|^2,  \qquad \psi\in
K^0_0,
\end{equation*}
where the norm $\|\cdot\|$ is given by expression (\ref{34}). As is
well known for such operators (see \cite[Theorem X.26]{RS}) the
essential self-adjointness of $H_{A,V}^c$ is equivalent to the
condition that $Ker(H_{min}^*)=\{0\}$. Here the kernel of this
operator is denoted by $Ker$. Then the essential self-adjointness of
$H_{A,V}^c$ means that the equation
\begin{equation}\label{40}
H_{A,V}^c\psi=0
 \end{equation}
 has only a
 trivial solution in $\mathcal{H}^0$.

 Let $\psi$ be a solution of equation
(\ref{40}). We introduce notation $\psi^N=\chi^N\cup\psi$ and
suppose that  $H_{A,V}^c\psi^N=f^N$. Then \begin{align*}
\big(H_{A,V}^c\psi^{\tau N}, \ \psi^{\tau
N}\big)_N&=\sum_{|k|,|s|\leq N}f_{k,s}^{\tau
N}\cdot\overline{\psi}^{\tau N}_{k,s}\geq\sum_{|k|,|s|\leq
N}|\psi_{k,s}^{\tau N}|^2=\\&=\sum_{|k|,|s|\leq
N}|\psi_{k,s}|^2=\|\psi^N\|^2.
\end{align*}
On the other hand, since $H_{A,V}^c\psi=0$ according to our
assumption, relation (\ref{39}) yields
 \begin{equation*}\big(H_{A,V}^c\psi^{\tau
N}, \ \psi^{\tau N}\big)_N=0.
\end{equation*}
Hence,
\begin{equation*}
\|\psi^N\|^2\leq 0.
\end{equation*}
  Passing to the limit as $N\rightarrow+\infty$, we obtain $\|\psi^N\|^2\rightarrow\|\psi\|^2$.
  
 Thus, $\psi=0$.
\end{proof}

\begin{cor} Suppose that the discrete electric potential $V\in K^0$ is
bounded from below, i.e., there exists $c$ such that, for all $k,
s\in \mathbb{Z}$, the inequality $V_{k,s}\geq c>-\infty$ is
satisfied. Then the operator $H_{A,V}^c$  is essentially
self-adjoint.
\end{cor}
\begin{proof} Indeed, since the discrete magnetic Laplacian
$-\Delta^c_A$ is a positive operator on $K^0_0$ (see proof in
\cite{Sushch}), the boundedness from below of the form $V\in K^0$
leads to the semiboundedness from below of the operator $H_{A,V}^c$.
\end{proof}

\end{document}